\documentclass[12pt,reqno]{amsart}
\usepackage{amsthm}
\usepackage{booktabs}

\usepackage{amssymb}
\usepackage{graphics}   
\usepackage{tikz}
\usetikzlibrary{shapes,backgrounds,calc}
\usepackage{latexsym}
\usepackage{multicol}
\usepackage{verbatim,enumerate}
\usepackage{accents}
\usepackage{cite}
\usepackage{enumitem}
\usepackage{array}
\usepackage{mathtools}
\usepackage{dsfont}
\usepackage{nicematrix}

\usepackage[colorlinks=true, linkcolor=blue, citecolor=blue, urlcolor=blue]{hyperref}
\usepackage{hyperref}
\usepackage{amsmath, amscd,url}
\usetikzlibrary{decorations.pathreplacing}

\usepackage{setspace}

\usepackage{pstricks}

\advance\textwidth by 1.5in 
\usepackage{geometry}
\newtheorem{thm}{Theorem}[section]
\newtheorem{cor}[thm]{Corollary}
\newtheorem{lem}[thm]{Lemma}
\newtheorem{prop}[thm]{Proposition}
\theoremstyle{definition}

\newtheorem{defn}[thm]{Definition}
\newtheorem{example}[thm]{Example}
\newtheorem{rem}[thm]{Remark}

\numberwithin{equation}{subsection}

\newenvironment{pf}{\proof}{\endproof}
\newcounter{cnt}
 \makeatletter
\def\mydggeometry{\makeatletter\dg@YGRID=1\dg@XGRID=20\unitlength=0.003pt\makeatother}
\makeatother \theoremstyle{remark}
\numberwithin{equation}{section}
\let\bwdg\bigwedge
\def\bigwedge{{\textstyle\bwdg}}

\newcommand{\thmref}[1]{Theorem~\ref{#1}}
\newcommand{\secref}[1]{Section~\ref{#1}}
\newcommand{\lemref}[1]{Lemma~\ref{#1}}
\newcommand{\propref}[1]{Proposition~\ref{#1}}
\newcommand{\corref}[1]{Corollary~\ref{#1}}
\newcommand{\remref}[1]{Remark~\ref{#1}}

\newcommand{\nc}{\newcommand}
\newcommand{\rnc}{\renewcommand}

\nc{\cal}{\mathcal} \nc{\goth}{\mathfrak} \rnc{\bold}{\mathbf}

\newcommand{\supp}{\operatorname{supp}}

\nc\bomega{{\mbox{\boldmath $\omega$}}} \nc\bpsi{{\mbox{\boldmath $\Psi$}}}
 \nc\balpha{{\mbox{\boldmath $\alpha$}}}
 \nc\bpi{{\mbox{\boldmath $\pi$}}}
 \nc\bvpi{{\mbox{\boldmath $\varpi$}}}
\nc\chara{\operatorname{ch}}

  \nc\bxi{{\mbox{\boldmath $\xi$}}}
\nc\bmu{{\mbox{\boldmath $\mu$}}} \nc\bcN{{\mbox{\boldmath $\cal{N}$}}} \nc\bcm{{\mbox{\boldmath $\cal{M}$}}} \nc\blambda{{\mbox{\boldmath
$\lambda$}}}\nc\bnu{{\mbox{\boldmath $\nu$}}}

\makeatletter
\def\section{\def\@secnumfont{\mdseries}\@startsection{section}{1}%
  \z@{.7\linespacing\@plus\linespacing}{.5\linespacing}%
  {\normalfont\scshape\centering}}
\def\subsection{\def\@secnumfont{\bfseries}\@startsection{subsection}{2}%
  {\parindent}{.5\linespacing\@plus.7\linespacing}{-.5em}%
  {\normalfont\bfseries}}
\makeatother

 \nc{\Hom}{\operatorname{Hom}}
  \nc{\mode}{\operatorname{mod}}
\nc{\End}{\operatorname{End}} \nc{\wh}[1]{\widehat{#1}} \nc{\Ext}{\operatorname{Ext}} \nc{\ch}{\text{ch}} \nc{\ev}{\operatorname{ev}}
\nc{\Ob}{\operatorname{Ob}} \nc{\soc}{\operatorname{soc}} \nc{\rad}{\operatorname{rad}} \nc{\head}{\operatorname{head}}

\def\det{\operatorname{det}}

 \nc{\Cal}{\cal} \nc{\Xp}[1]{X^+(#1)} \nc{\Xm}[1]{X^-(#1)}
\nc{\on}{\operatorname} \nc{\Z}{{\bold Z}} \nc{\J}{{\cal J}}  \nc{\Q}{{\bold Q}}

\nc{\N}{{\bold N}}  \nc\boa{\bold a} \nc\bob{\bold b} \nc\boc{\bold c} \nc\bod{\bold d} \nc\boe{\bold e} \nc\bof{\bold f} \nc\bog{\bold g}
\nc\boh{\bold h} \nc\boi{\bold i} \nc\boj{\bold j} \nc\bok{\bold k} \nc\bol{\bold l} \nc\bom{\bold m} \nc\bon{\mathbb n} \nc\boo{\bold o}
\nc\bop{\bold p} \nc\boq{\bold q} \nc\bor{\bold r} \nc\bos{\bold s} \nc\boT{\bold t} \nc\boF{\bold F} \nc\bou{\bold u} \nc\bov{\bold v}
\nc\bow{\bold w} \nc\boz{\bold z}\nc\ba{\bold A} \nc\bb{\bold B} \nc\bc{\mathbb C} \nc\bd{\bold D} \nc\be{\bold E} \nc\bg{\bold
G} \nc\bh{\bold H} \nc\bi{\bold I} \nc\bj{\bold J} \nc\bk{\bold K} \nc\bl{\bold L} \nc\bm{\bold M} \nc\bn{\mathbb N} \nc\bo{\bold O} \nc\bp{\bold
P} \nc\bq{\bold Q} \nc\br{\bold R} \nc\bs{\bold S} \nc\bt{\bold T} \nc\bu{\bold U} \nc\bv{\bold V} \nc\bw{\bold W} \nc\bz{\mathbb Z} \nc\bx{\bold
x} \nc\KR{\bold{KR}} \nc\rk{\bold{rk}} \nc\het{\text{ht }}

\nc\toa{\tilde a} \nc\tob{\tilde b} \nc\toc{\tilde c} \nc\tod{\tilde d} \nc\toe{\tilde e} \nc\tof{\tilde f} \nc\tog{\tilde g} \nc\toh{\tilde h}
\nc\toi{\tilde i} \nc\toj{\tilde j} \nc\tok{\tilde k} \nc\tol{\tilde l} \nc\tom{\tilde m} \nc\ton{\tilde n} \nc\too{\tilde o} \nc\toq{\tilde q}
\nc\tor{\tilde r} \nc\tos{\tilde s} \nc\toT{\tilde t} \nc\tou{\tilde u} \nc\tov{\tilde v} \nc\tow{\tilde w} \nc\toz{\tilde z} \nc\woi{w_{\omega_i}}

\newcommand\wrapped[1]%
{\renewcommand\arraystretch{1}%
	\begin{array}{@{}l@{}}#1\end{array}%
}

\ifTUTeX
  \usepackage{fontspec}
\else
  \usepackage[T1]{fontenc}
  \usepackage[utf8]{inputenc} % The default since 2018
  \DeclareUnicodeCharacter{200B}{{\hskip 0pt}}
\fi

\begin{document}
	
	\title{The Full P-vertex Problem and Perfect Matchings for Bipartite Graphs}

   \author{G. Arunkumar}
	
	\address{Indian Institute of Technology Madras, Chennai, India.}
	\email{garunkumar@iitm.ac.in}
 
	\author{Puja Samanta}
    \thanks{}
	\address{Indian Institute of Technology Madras, Chennai, India.}
	\email{ma23d004@smail.iitm.ac.in, pujasamanta1999@gmail.com}

	%\thanks{}

	%\thanks{}

	%\thanks{}
	%\thanks{$^{\ast}$ supported by the fellowship 0203/8(28)/2023-R\&D-II from the National Board for Higher Mathematics (NBHM), Govt. of India. }
   %\thanks{$^{1}$ supported by the NFIG grant of the Indian Institute of Technology Madras RF/22-23/0985/MA/NFIG/009003 and the SERB Startup Research Grant SRG/2022/001281.} 


\begin{abstract}
In a recent work, Sharma and Panda~\cite{sharma} showed that every bipartite graph with a perfect matching has property (P), and proved the converse for trees and unicyclic bipartite graphs (i.e., bipartite graphs with cycle rank $m(G) \le 1$). In this paper, we extend this result to broader classes of bipartite graphs. We first show that every bipartite graph with property (P) is balanced. We prove that the converse holds for all bipartite graphs with cycle rank $m(G)$ at most three, and further establish it for several additional families of bipartite graphs. Finally, we derive algebraic constraints for balanced bipartite graphs without perfect matchings and use them to identify a family of bipartite graphs that does not have property (P).
\end{abstract}

\maketitle

\medskip

\noindent\textbf{Keywords:}
Property (P); P-vertex; biadjacency matrix; perfect matching; cycle rank; hollow matrix.
%\maketitle
%\tableofcontents
%\clearpage

\section{Introduction}
Let $G$ be a simple undirected graph with vertex set $V(G)$ and edge set $E(G)$. Throughout this paper, all  matrices are assumed to be real. For a matrix $A=(A_{ij})\in\mathbb{R}^{n\times n}$, we write $A^\top$ for its transpose. If $\lambda$ is an eigenvalue of $A$, then $m_A(\lambda)$ denotes its algebraic multiplicity. 
The support graph of a symmetric matrix $A$, denoted by $G(A)$, is the graph on the vertex set $[n]=\{1,\dots,n\}$ with edge set
\[
E(G(A))=\bigl\{ij : i\neq j\ \text{and}\ A_{ij}\neq 0\bigr\}.
\]
Given a graph $G$ on $n$ vertices, define
\[
S(G):=\{A\in\mathbb{R}^{n\times n}: A^\top=A\ \text{and}\ G(A)=G\}.
\]
A matrix $A \in S(G)$ is called an (irreducible) acyclic matrix
if the graph $G$ is a tree.
For a nonempty index set $\alpha\subseteq[n]$, let $A(\alpha)$ denote the principal submatrix of $A$ obtained by deleting the rows and columns indexed by $\alpha$ (for $\alpha=\{i\}$ we write $A(i)$).
By Cauchy interlacing theorem \cite[Theorem 4.3.17]{inter}, for each $i\in[n]$ and eigenvalue $\lambda$ of $A$,
\[
m_{A(i)}(\lambda)-m_A(\lambda)\in\{-1,0,1\}.
\]
It is convenient to write
\[
\epsilon_i(\lambda):=m_{A(i)}(\lambda)-m_A(\lambda).
\]
Recall that a vertex $i \in V(G)$ is called a Fiedler-vertex of $A$ for the eigenvalue $\lambda$ if 
$\epsilon_i(\lambda) \in \{0,1\}$ \cite{Fiedler,hermitian}. 
When $\lambda = 0$, such a vertex is called an F-vertex. 
Moreover, if $\epsilon_i(\lambda) = 1$, then $i$ is called a Parter-vertex of $A$ for $\lambda$; 
when $\lambda = 0$, it is called a P-vertex \cite{kim1}. An F-vertex that is not a P-vertex (i.e., $\epsilon_i(0)=0$) is called a neutral vertex, and a vertex with $\epsilon_i(\lambda)=-1$ is called a downer vertex \cite{hermitian}. The total number of P-vertices of $A$ is denoted by $P_\nu(A)$.
A graph $G$ of order $n$ is said to have full P-vertices if there exists $A\in S(G)$ with $P_\nu(A)=n$.

More generally, for any nonempty subset $\alpha\subseteq[n]$, by the Cauchy interlacing theorem, we have
\(
m_{A(\alpha)}(0)\le m_A(0)+|\alpha|.
\)
If the equality holds, then $\alpha$ is called a P-set \cite{kim2}. 

The study of P-vertices can be traced back to the foundational work of Parter \cite{Parter}, who investigated how the structure of a graph, especially a tree, affects the multiplicities of eigenvalues of associated real symmetric matrices. Wiener \cite{wiener} later expanded and enriched this line of research by providing a comprehensive analysis of the location and multiplicity of eigenvalues of sign symmetric acyclic matrices. The formal introduction and systematic study of P-vertices appeared subsequently in \cite{hermitian, kim1}, where the concept was refined and explored in greater depth. Godsil has also studied an analogous notion within the framework of the matching polynomial theory, where such a vertex appears as a $\theta$-positive vertex of a graph $G$ \cite{godsil}. 
%%%%%%%%%%%%%%%%%%%%%%%%%%%%%%%%%%%%%%%%%%%%%%%%%%%%%%%%%%%%%%%%%%%%%%%%%%%
Johnson and Sutton \cite{hermitian} initiated a systematic investigation of how vertex deletion affects eigenvalue multiplicities of Hermitian matrices. Furthermore, they showed that every singular acyclic matrix of order $n$ has at most $n-2$ number of P-vertices. Kim and Shader \cite{kim1} later proved that this bound need not extend to nonsingular acyclic matrices by exhibiting a counterexample. Moreover, for a nonsingular tridiagonal matrix $A$ of order $n$ they \cite[Theorem 4]{kim1} established
\[
P_\nu(A)\le
\begin{cases}
n-1, & \text{if $n$ is odd},\\[2pt]
n,   & \text{if $n$ is even}.
\end{cases}
\]
In particular, for such a matrix $A$, the equality $P_{\nu}(A) = n$ is possible only when $n$ is even.

These results motivated a systematic study of P-vertices for acyclic graphs. Several open problems posed by Kim and Shader \cite{kim1} were addressed and resolved in \cite{fons1,fons2,fons3,fons4}. Trees that admit a nonsingular matrix in which every vertex is a P-vertex were structurally characterized in \cite{fons2,fons5}.

In 2021, Fonseca et al.~\cite{fons6} showed that the number of P-vertices of a singular acyclic matrix $A$ of rank $r$ is at most $r-2$ when $r$ is odd. Subsequently, Du and Fonseca~\cite{du2} characterized all trees on $n$ vertices that admit a singular matrix with nullity $k$, proving that the number of P-vertices is either $n-k-1$ or $n-k-2$. Further results on bounds for the number of P-vertices in trees were obtained in \cite{cruz,du5}. These results motivate further study of P-vertices. Moreover, in this direction, the notion of P-sets has been studied in the literature (see, e.g., \cite{du1,du3,du4}).

We now turn to the problem of determining when a graph admits a matrix in which every vertex is a P-vertex, particularly beyond the class of trees. This problem has been studied for several graph families. It was shown in \cite[Theorem 4.2]{fons1} that every cycle graph $C_n$ has the full P-vertex property. In 2025, Howlader et al.\ \cite[Definition 1.1]{Howlader} introduced the terminology property (P) as a reformulation of the full P-vertex property, namely, a graph $G$ has property (P) if there exists a nonsingular
symmetric matrix $A \in S(G)$ for which every vertex is a P-vertex. They further extended this line of investigation to unicyclic graphs \cite[Theorem~6.1]{Howlader}. 
More recently, in 2026, Sharma and Panda ~\cite{sharma} studied property (P) for bipartite graphs and established the following result for connected bipartite graphs.
 \begin{thm}\cite[Theorem 3.1]{sharma} \label{bipartite perfect matching theorem}
Let $G$ be a connected bipartite graph with a perfect matching. Then \(G\) has property (P).
\end{thm}
By a result of Howlader et al.\ \cite[Lemma~2.1]{Howlader}, a graph has property (P) if and only if each of its connected components has property (P). Likewise, a graph has a perfect matching if and only if each of its connected components has a perfect matching. Therefore, Theorem~\ref{bipartite perfect matching theorem} naturally extends to arbitrary (not necessarily connected) bipartite graphs. Henceforth, we apply Theorem~\ref{bipartite perfect matching theorem} directly to arbitrary bipartite graphs, with the connected case understood componentwise.

Moreover, for acyclic and unicyclic bipartite graphs, which correspond to cycle rank $m(G)\le 1$, the existence of a perfect matching is both necessary and sufficient for property (P)~\cite[Theorems~3.3 and~3.8]{sharma}. However, in general, it remains open whether every bipartite graph with property (P) has a perfect matching. The converse has been established only under additional assumptions~\cite[Theorem~3.12]{sharma}.

In this paper, we study property (P) for bipartite graphs with the aim of understanding when it forces the existence of a perfect matching. We first show that every bipartite graph with property (P) must be balanced (\thmref{thm:balanced}), which immediately excludes bipartite graphs of odd order. As a direct application of this structural requirement, we completely determine when complete bipartite graphs and grid graphs possess property (P). We next extend the known converse for acyclic and unicyclic bipartite graphs to bipartite graphs of cycle rank at most three, thereby including bicyclic and tricyclic graphs (\thmref{thm:cyclerank3}). We also identify a graph of cycle rank four that motivates the algebraic techniques developed later in the paper. 
Furthermore, we prove that property (P) is equivalent to the existence of a perfect matching for several additional classes of bipartite graphs, including bipartite graphs of order $8+2k \ (k \ge 0)$ with at least $k$ pendant vertices (\thmref{8+2k}), bipartite graphs that can be reduced to order at most $8$ under the pendant-vertex reduction process (\thmref{pendant--edge reduction theorem}), and triangular bipartite graphs (\thmref{triangular-theorem}).
Together, these results describe several natural families of bipartite graphs
for which property (P) is equivalent to the existence of a perfect matching. 

To further investigate the remaining open cases of the converse, we consider balanced bipartite graphs that do not admit a perfect matching in \secref{sec:matrix_constraints}. Using a block decomposition arising from a violation of Hall's condition, we derive necessary algebraic constraints that every hollow-inverse matrix on such a graph must satisfy. These constraints yield a structural obstruction for a family of bipartite graphs, showing that none of them has property (P) (\thmref{thm:generalized_obstruction}). We also prove that any such matrix associated with a bipartite graph lacking a perfect matching must have at least two nonzero diagonal entries (\propref{lem:nondeg}), complementing the result of Sharma and Panda~\cite[Theorem~3.12]{sharma} concerning hollow-inverse matrices with at most three nonzero diagonal entries. 

Throughout the paper, we follow the graph-theoretic terminology of \cite{west}.

%%%%%%%%%%%%%%%%%%%%%%%%%%%%%%%%%%%%%%%%%%%%%%%%%%%%%%%%%%%%%%%%%%%%%%%%%%%%%%%%%%%%%%%%%%%%%%%%%%%%%%%%%%%%%%%%%%%%%%%%%%

\section{Necessary Condition for Property~(P)}\label{sec:converse section}
In this section, we establish the first structural consequence of property (P), showing that every bipartite graph with property (P) must be balanced.

We adopt the following terminology.
\begin{defn}
A nonsingular matrix $A$ is called a \emph{hollow-inverse matrix} if all the diagonal entries of $A^{-1}$ are zero.
\end{defn}
In particular, if $G$ is a graph and $A\in S(G)$, then $A$ is a hollow-inverse matrix if and only if every vertex of $G$ is a P-vertex of $A$. Consequently, a graph has property~(P) if and only if it admits a hollow-inverse matrix.

We denote by $G=G(U,V;E)$ a bipartite graph with partite sets \( U = \{u_1,\dots,u_m\},\) and  \( V= \{v_1,\dots,v_n\}\). We say that $G$ is a balanced bipartite graph if $|U|=|V|$, and unbalanced otherwise.
The biadjacency matrix $B=(B_{ij})\in\mathbb R^{m\times n}$ of $G$ is defined by
\[
B_{ij}=
\begin{cases}
1,&\text{if }u_iv_j\in E,\\
0,&\text{otherwise.}
\end{cases}
\]
For any matrix $M=(M_{ij}) \in \mathbb{R}^{m \times n}$, the \emph{support} of $M$ is 
    \[
    \supp(M) =\{(i,j)\in [m] \times [n] : M_{ij}\neq 0\}.
    \]
Since an $m\times n$ matrix naturally defines a bipartite graph whose edges correspond to its nonzero entries \cite[p.~3]{du4}, we make the following definition.
\begin{defn}
A matrix $M\in\mathbb{R}^{m\times n}$ has the
\emph{zero--nonzero pattern defined by $G$} if
\[
\supp(M)=\{(i,j):u_iv_j\in E\}.
\]
Equivalently,
\[
\supp(M)=\supp(B).
\]
\end{defn}
By ordering the vertices of $U$ before those of $V$, any symmetric matrix $A\in S(G)$ naturally takes the block form
\[
\renewcommand{\arraystretch}{1.3}
A=
\left[
\begin{array}{c|c}
X & D\\ \hline
D^{\top} & Y
\end{array}
\right],
\]
where $X\in\mathbb{R}^{m\times m}$ and $Y\in\mathbb{R}^{n\times n}$ are diagonal matrices, and $D\in\mathbb{R}^{m\times n}$ has the zero-nonzero pattern defined by $G$.
%----------------------------------------------------------------------------------------------
\begin{thm}\label{thm:balanced}
    Let $G$ be a bipartite graph. If $G$ has property (P), then $G$ is a balanced bipartite graph.
\end{thm}
\begin{pf}
    Let the bipartition of $G$ be $U=\{u_1,u_2,\dots,u_m\}$ and $V=\{v_1,v_2,\dots,v_n\}$. 
    
    Since $G$ has property (P), it admits a hollow-inverse matrix $A \in S(G)$. Using the notation introduced above, write
\[
\renewcommand{\arraystretch}{1.3}
A=
\left[
\begin{array}{c|c}
X & D\\ \hline
D^{\top} & Y
\end{array}
\right]
 \ \text{ and }  \ 
\renewcommand{\arraystretch}{1.3}
A^{-1} =
\left[
\begin{array}{c|c}
P & Q\\ \hline
Q^{\top} & R
\end{array}
\right],
\]
where $P\in\mathbb{R}^{m\times m}$, $R\in\mathbb{R}^{n\times n}$ and $Q \in \mathbb{R}^{m\times n}$. 
Because $A$ is a hollow-inverse matrix, $A^{-1}$ has zero diagonal ($\operatorname{diag}(A^{-1})=\mathbf{0}$). Consequently, the diagonal entries of both principal blocks $P$ and $R$ are zero.

From $AA^{-1}=I_{m+n}$, we obtain the following block equations,
\begin{equation}\label{m diag}
    XP + DQ^{\top} = I_m,
\end{equation}
\begin{equation}\label{n diag}
    D^{\top}Q +YR = I_n,
\end{equation}
\begin{equation*}
    XQ + DR = \mathbf{0}_{m\times n},
\end{equation*}
\begin{equation*}
    D^{\top}P + YQ^{\top}= \mathbf{0}_{n \times m},
\end{equation*}
 where
$\mathbf{0}_{m \times n}$ denotes the zero matrix of order $m \times n$.

Taking the trace of \eqref{m diag} and using $(XP)_{ii} = 0$ for all
$i \in \{1,2,\dots,m\}$ (since $X$ is a diagonal matrix and $P$ has zero diagonal entries),
we obtain
\begin{equation}\label{m trace}
    \operatorname{Tr}(XP + DQ^{\top})= \operatorname{Tr} (DQ^{\top}) = m.
\end{equation}
Similarly, taking the trace of \eqref{n diag} and using $(YR)_{jj} = 0$ for all
$j \in \{1,2,\dots,n\}$ (since $Y$ is a diagonal matrix and $R$ has zero diagonal entries),
we obtain
\begin{equation}\label{n trace}
    \operatorname{Tr}(D^{\top}Q +YR) = \operatorname{Tr}(D^{\top}Q)= n.
\end{equation}
From \eqref{m trace} and \eqref{n trace},
\[
m = \operatorname{Tr}(DQ^{\top}) = \operatorname{Tr}(QD^{\top})
= \operatorname{Tr}(D^{\top}Q) = n.
\]
Hence, $m=n$, so $|U|=|V|$. Therefore, $G$ is a balanced bipartite graph.
\end{pf}
As an immediate consequence, we obtain the following.
\begin{cor}\label{cor:odd bipartite}
No bipartite graph of odd order has property (P).
\end{cor}
We next record two direct applications of the preceding results to complete bipartite graphs and grid graphs.

\begin{cor}\label{cor:complete bipartite graph}
The complete bipartite graph $K_{m,n}$ has property (P) if and only if $m=n$.
\end{cor}
\begin{pf}
If $m=n$, then the complete bipartite graph $K_{n,n}$ has a perfect matching, and hence by \thmref{bipartite perfect matching theorem} it has property (P). Conversely, if the graph $K_{m,n}$ has property (P), then by \thmref{thm:balanced} it must be balanced.
\end{pf}
\begin{rem}
We give an explicit matrix realizing property (P) for the graph $K_{n,n}$. Let
\[
U=\{u_1,\dots,u_n\}, \
V=\{v_1,\dots,v_n\},
\]
and fix the perfect matching
\[
M=\{u_i v_i :1\le i\le n\}.
\]
Define the matrix
\[
D=(1-w)I_n+wJ_n,
\]
where $J_n$ denotes the $n\times n$ all-ones matrix and $w\neq1,-\frac1{n-1}$. That is,
\[
D_{ij}=
\begin{cases}
1,& i=j,\\
w,& i\neq j.
\end{cases}
\]
Since
\[
\det(D)=(1-w)^{n-1}\bigl(1+(n-1)w\bigr)\neq0,
\]
the matrix $D$ is invertible. Define
\[
\renewcommand{\arraystretch}{1.3}
A=
\left[
\begin{array}{c|c}
\mathbf{0}_n & D\\ \hline
D^{\top} & \mathbf{0}_n
\end{array}
\right],
\]
where $\mathbf{0}_n$ denotes the zero matrix of order $n$.
Then $A\in S(K_{n,n})$ is nonsingular and
\[
\renewcommand{\arraystretch}{1.3}
A^{-1}=
\left[
\begin{array}{c|c}
\mathbf{0}_n & (D^{\top})^{-1}\\ \hline
D^{-1} & \mathbf{0}_n
\end{array}
\right].
\]
Since the diagonal entries of $A^{-1}$ are all zero, $A$ is a hollow-inverse matrix for $K_{n,n}$.
\end{rem}
We next consider another standard family of bipartite graphs and obtain a complete characterization of property (P).

An $m\times n$ grid graph, denoted by $G_{m,n}$, has vertex set
\[
V=\{(i,j) : 1\le i\le m,\ 1\le j\le n\},
\]
arranged in $m$ rows and $n$ columns.
Two vertices $(i,j)$ and $(i',j')$ are adjacent
whenever $|i-i'|+|j-j'|=1$.   For example, the grid graph $G_{3,4}$ is illustrated in Figure \ref{fig:grid23}. 
We recall the following standard characterization of perfect matchings in grid graphs \cite[p.~48]{grid perfect matching reference}.

\begin{prop}\label{prop:gridmatching}
The grid graph $G_{m,n}$ is bipartite, and it has a perfect matching if and only if $mn$ is even.
\end{prop}
Combining \propref{prop:gridmatching} with our earlier results yields the following characterization.
\begin{cor}
The grid graph $G_{m,n}$ has property (P) if and only if $mn$ is even.
\end{cor}
\begin{pf}
By \propref{prop:gridmatching}, the graph $G_{m,n}$ has a perfect matching if and only if $mn$ is even. The result now follows from \thmref{bipartite perfect matching theorem} and \corref{cor:odd bipartite}.
\end{pf}
\begin{figure}[h!]
\centering
\begin{tikzpicture}[
  vertex/.style={circle, draw=black, fill=black, inner sep=1.6pt},
  edge/.style={gray!60},
  matching/.style={line width=1.4pt, blue},
  scale=1.3
]

\def\cols{4}
\def\rows{3}

\foreach \i in {1,...,\cols}{
  \foreach \j in {1,...,\rows}{
    \coordinate (v\i\j) at (\i,\j);
    \node[vertex] (n\i\j) at (v\i\j) {};
  }
}

\foreach \i in {1,...,\numexpr\cols-1\relax}{
  \foreach \j in {1,...,\rows}{
    \draw[edge] (n\i\j)--(n\the\numexpr\i+1\relax\j);
  }
}

\foreach \i in {1,...,\cols}{
  \foreach \j in {1,...,\numexpr\rows-1\relax}{
    \draw[edge] (n\i\j)--(n\i\the\numexpr\j+1\relax);
  }
}

\draw[matching] (n11)--(n21);
\draw[matching] (n31)--(n41);

\draw[matching] (n12)--(n13);
\draw[matching] (n22)--(n23);
\draw[matching] (n32)--(n33);
\draw[matching] (n42)--(n43);

\end{tikzpicture}
\caption{The grid graph $G_{3,4}$ with $3\times4$ vertices
and admits a perfect matching, shown in blue.}
\label{fig:grid23}
\end{figure}
\section{Bipartite Graphs of Low Cycle Rank}\label{sec:schur_reductions}
We now relate property (P) to the cycle rank of bipartite graphs. After introducing the notion of cycle rank, we prove that property (P) is equivalent to the existence of a perfect matching for every bipartite graph of cycle rank at most three.

Recall that the \emph{edge space} of a graph \(G\), denoted by \(\mathcal{E}(G)\), is the vector space $\{0,1\}^E$  over \(\mathbb{F}_2\), which we view as the power set of the edge set \(E\), with symmetric difference as addition. 
% Thus, $\Big\{ \{e_1\},\dots,\{e_{|E|}\} \Big\}$ is the standard basis of \(\mathcal{E}(G)\) and 
% \(
% \dim \mathcal{E}(G)=|E|.
% \)
The \emph{cycle space} of \(G\), denoted by \(\mathcal{C}(G)\), is the subspace of \(\mathcal{E}(G)\) spanned by all the cycles in $G$, more precisely by the edge sets of all cycles of \(G\) \cite[p.~23-24]{diestel}.
\begin{defn}[{\cite[p.~24]{diestel}, \cite[p.~39]{harary}}]\label{defn:CR}
The dimension of the cycle space $\mathcal{C}(G)$ is called the \emph{cycle rank} (or \emph{cyclomatic number}) of $G$, denoted by $m(G)$.
\end{defn}
% Following Harary~\cite[p.~38]{harary}, a collection of cycles $\{Z_1, Z_2, \dots, Z_k\}$ in $G$ is said to be \emph{linearly independent} if no cycle in the collection can be expressed as a linear combination (under symmetric difference over $\mathbb{F}_2$) of the remaining cycles. A \emph{cycle basis} of $G$ is a maximal collection of linearly independent cycles.
The following fundamental formula expresses the cycle rank in terms of the numbers of vertices, edges, and connected components of a graph. 
\begin{thm}\cite[Theorem~4.5(b)]{harary}\label{thm:CR}
  For every graph $G$,
\[
m(G)=|E(G)|-|V(G)|+k(G),
\]
where $k(G)$ denotes the number of connected components of $G$.  
\end{thm}
\begin{rem}\label{rem:graph_classes}
As an immediate consequence of Theorem~\ref{thm:CR}, a graph is acyclic if and only if $m(G)=0$. Since a graph is unicyclic if and only if it contains exactly one cycle, it follows that $m(G)=1$ in this case. 
Following standard terminology (see, e.g., Cvetkovi\'{c} et al.~\cite[p.~6]{cvetkovic2010}), a connected graph $G$ with $|V(G)|+1$ edges is called \emph{bicyclic}. Analogously, we refer to a connected graph with $|V(G)|+2$ edges as \emph{tricyclic}. By Theorem~\ref{thm:CR}, these are precisely the connected graphs with cycle ranks $2$ and $3$, respectively.
\end{rem}

For a bipartite graph $G(U,V;E)$ and a subset $S\subseteq U$, let
\[
N(S):=\{v\in V:\exists\,u\in S \text{ such that } uv \in E\}
\]
denote the set of neighbors of $S$. For subsets $X,Y\subseteq V(G)$, let $e(X,Y)$ denote the number of edges with one endpoint in $X$ and the other in $Y$. We now recall Hall's theorem, which characterizes the existence of perfect matchings in bipartite graphs.
\begin{thm}[Hall's Theorem, {\cite[Theorem~3.1.11]{west}} ]\label{hall}
    A bipartite graph $G(U,V;E)$, with $|U|=|V|$, has a perfect matching if and only if
    \[
    |N(S)|\geq |S|,
    \]
    for all $S \subseteq U$.
\end{thm}
For a vertex \(v\in V(G)\), let \(G-v\) denote the graph obtained by deleting the vertex \(v\) and all the edges incident to it; consequently, $G - \{u,v\}$ denotes the sequential deletion of both the vertices $u$ and $v$. The following result of Howlader et al.~\cite[Theorem 2.3]{Howlader} will be useful. For consistency with our framework, we restate the result using the deletion notation.
\begin{thm}\cite[Theorem 2.3]{Howlader}\label{lem:pendant_heredity}
Let $G$ be a graph. If $u$ is a pendant vertex of $G$ with the unique neighbor $v$, then $G - \{u,v\}$ has property (P) if and only if $G$ has property (P).
\end{thm}
To prove the characterization for bipartite graphs of cycle rank at most three, we first establish two lemmas.
%%%%%%%%%%%
\begin{lem}\label{lem:isolated_vertex}
If a graph $G$ contains an isolated vertex, then $G$ does not have property (P). Consequently, any graph possessing property (P) must satisfy $\delta(G) \ge 1$.
\end{lem}
\begin{pf}
    Suppose, to the contrary, that $G$ has property (P). Then there exists a nonsingular matrix $A\in S(G)$ such that $A^{-1}$ is hollow. Let $v$ be an isolated vertex of $G$, and let $e_v$ denote the corresponding standard basis vector, then 
    \(
Ae_v=A_{vv}e_v.
\)
As $A$ is nonsingular, we have $A_{vv}\neq 0$. Hence
\[
A^{-1}e_v=A_{vv}^{-1}e_v,
\]
which implies
\(
(A^{-1})_{vv}=A_{vv}^{-1}\neq0.
\) 
This contradicts the assumption that $A^{-1}$ is hollow. Therefore, $G$ cannot have property (P).
\end{pf}
\begin{lem}\label{prop:size}
Let $H$ be a connected balanced bipartite graph with minimum degree
$\delta(H)\ge 2$. If $H$ has no perfect matching, then
\(
    m(H) \ge 4.
\)
\end{lem}

\begin{pf}
Let $(U,V)$ be the bipartition of $H$, with $|U|=|V|=n$. Since $H$ has no
perfect matching, Hall's theorem gives a subset $U_1\subseteq U$ such that
\[
    |N(U_1)|<|U_1|.
\]
Set
\[
    V_1 =N(U_1), \qquad U_2 =U\setminus U_1, \qquad V_2 =V\setminus V_1.
\]
Let $|U_1|=r$ and $|V_1|=k$. Then $r>k$. Put $d=r-k\ge 1$.

Every edge incident with a vertex of $U_1$ has its other endpoint in $V_1$.
Hence, using $\delta(H)\ge 2$, we get
\[
    e(U_1,V_1)\ge 2r.
\]
Moreover, since $V_1=N(U_1)$, there are no edges between $U_1$ and $V_2$.
Thus, every edge incident with a vertex of $V_2$ has its other endpoint in
$U_2$. Since $|V_2|=n-k$ and $\delta(H)\ge 2$, we obtain
\[
    e(U_2,V_2)\ge 2(n-k).
\]
Also, $V_2\neq \emptyset$, because $k<r\le n$.

Finally, since $H$ is connected, there must be at least one edge connecting
$U_1\cup V_1$ to $U_2\cup V_2$. This edge cannot be between $U_1$ and $V_2$,
and hence must be between $U_2$ and $V_1$. Therefore,
\[
    e(U_2,V_1)\ge 1.
\]
The edges counted by $e(U_1,V_1)$, $e(U_2,V_2)$, and $e(U_2,V_1)$ are pairwise disjoint. Consequently,
\[
\begin{aligned}
    |E(H)|
    &\ge e(U_1,V_1)+e(U_2,V_2)+e(U_2,V_1) \\
    &\ge 2r+2(n-k)+1 \\
    &= 2n+2(r-k)+1 \\
    &= 2n+2d+1 \\
    &\ge 2n+3.
\end{aligned}
\]
Since $|V(H)|=2n$, this proves that
\[
    |E(H)|\ge |V(H)|+3.
\]
By \thmref{thm:CR}, and since $H$ is connected,
\[
m(H)=|E(H)|-|V(H)|+1\ge 4.
\]
\end{pf}
%%%%%%%%%%%
%%%%%%%%%%%%%%%%%%%
\begin{thm}\label{thm:cyclerank3}
Let $G$ be a bipartite graph with $m(G)\le 3$. Then $G$ has
property~(P) if and only if $G$ has a perfect matching.
\end{thm}

\begin{pf}
The sufficiency follows from \thmref{bipartite perfect matching theorem}.

Conversely, suppose that $G$ has property~(P). Then there exists a
hollow-inverse matrix $A\in S(G)$. We prove that $G$ has a perfect matching
by induction on $|V(G)|$.

By Lemma~\ref{lem:isolated_vertex}, $G$ has no isolated vertices. Hence
$\delta(G)\ge1$.

The result is immediate when $|V(G)|=2$.

Suppose that $G$ has a pendant vertex $u$, and let $v$ be its unique
neighbor. By \thmref{lem:pendant_heredity}, the graph $G-\{u,v\}$ also has
property~(P). Moreover, deleting vertices cannot increase the cycle rank, so
\[
    m(G-\{u,v\})\le m(G)\le 3.
\]
By the induction hypothesis, $G-\{u,v\}$ has a perfect matching $M'$. Hence
\[
    M'\cup \{uv\}
\]
is a perfect matching of $G$.

We may therefore assume that
\(
    \delta(G)\ge 2.
\)
Write the connected components of $G$ as
\(
G_1,G_2,\dots,G_t.
\)
By \cite[Lemma~2.1]{Howlader}, every connected component $G_i$ has property~(P). Hence, by \thmref{thm:balanced}, every $G_i$ is balanced.

Suppose, for a contradiction, that some component $H$ of $G$ has no perfect
matching. Then $H$ is a connected balanced bipartite graph with
\(
    \delta(H)\ge 2
\)
and no perfect matching. Hence, by \lemref{prop:size},
\[
    m(H)\ge 4.
\]
Since the cycle rank is nonnegative and additive over connected components,
\[
    m(H)\le m(G)\le 3,
\]
which is a contradiction.

Hence, every connected component of $G$ has a perfect matching. The union of these perfect matchings is a perfect matching of $G$. 
\end{pf}
%%%%%%%%%%%%%%%%%%%
\begin{rem}\label{rem:sharp}
 The graph $G$ in Figure~\ref{fig:sharp-example}
attains equality in \lemref{prop:size}, having $10$ vertices,
$13$ edges and $m(G)=4$. Consequently, this graph represents the
smallest possible obstruction beyond the range of graphs covered by
Theorem~\ref{thm:cyclerank3}. We return to this graph in Subsection~\ref{specific_example}, where we apply matrix-theoretic tools to prove that it does not have property~(P) (\corref{cor:rank_contradiction}).
\end{rem}

\begin{figure}[h]
\centering
\begin{tikzpicture}[scale=1]

% vertex style
\tikzset{
  vtx/.style={circle, fill=black, inner sep=1.6pt}
}

%---------------- Left part U ----------------%
\node[vtx] (u1) at (0,4) {}; \node[left=6pt] at (u1) {$u_1$};
\node[vtx] (u2) at (0,3) {}; \node[left=6pt] at (u2) {$u_2$};
\node[vtx] (u3) at (0,2) {}; \node[left=6pt] at (u3) {$u_3$};
\node[vtx] (u4) at (0,1) {}; \node[left=6pt] at (u4) {$u_4$};
\node[vtx] (u5) at (0,0) {}; \node[left=6pt] at (u5) {$u_5$};

%---------------- Right part W ----------------%
\node[vtx] (v1) at (3,4) {}; \node[right=6pt] at (v1) {$v_1$};
\node[vtx] (v2) at (3,3) {}; \node[right=6pt] at (v2) {$v_2$};
\node[vtx] (v3) at (3,2) {}; \node[right=6pt] at (v3) {$v_3$};
\node[vtx] (v4) at (3,1) {}; \node[right=6pt] at (v4) {$v_4$};
\node[vtx] (v5) at (3,0) {}; \node[right=6pt] at (v5) {$v_5$};

%---------------- Edges ----------------%

% K_{3,2}
\draw (u1)--(v1);
\draw (u1)--(v2);

\draw (u2)--(v1);
\draw (u2)--(v2);

\draw (u3)--(v1);
\draw (u3)--(v2);

% K_{2,3}
\draw (u4)--(v3);
\draw (u4)--(v4);
\draw (u4)--(v5);

\draw (u5)--(v3);
\draw (u5)--(v4);
\draw (u5)--(v5);

% bridge
\draw (u4)--(v1);

\end{tikzpicture}
\caption{A connected balanced bipartite graph $G$ of order $10$, with $\delta (G) = 2$ and $m(G)=4$.}
\label{fig:sharp-example}
\end{figure}

\section{Extending the Characterization of Property~(P) via Pendant-Vertex Reduction}\label{sec:small order}
We first establish the characterization of property (P) for bipartite graphs of small order and then extend it to a broader family of bipartite graphs containing pendant vertices. The following result from \cite{Howlader} will be useful in the study of bipartite graphs containing pendant vertices.
\begin{thm}\cite[Theorem~4.3]{Howlader}\label{antenna}
An antenna graph is a graph containing a vertex adjacent to more than one pendant vertex. No antenna graph has property (P).
\end{thm}
\begin{lem}\label{thm:order8}
   Let $G$ be a connected bipartite graph of order at most $8$. Then $G$ has property
(P) if and only if $G$ has a perfect matching.
\end{lem}
\begin{pf}
Let $G=G(U,V;E)$ be a connected bipartite graph. If $G$ admits a perfect
matching, then \thmref{bipartite perfect matching theorem} implies that $G$
has property~(P).

Conversely, suppose that $G$ has property~(P). By
\thmref{thm:balanced}, $G$ is balanced. Suppose, for a contradiction, that
$G$ does not have a perfect matching. By Hall's theorem, there exists a
subset $S\subseteq U$ such that
\[
|N(S)|<|S|.
\]

Since $G$ is balanced and has order at most $8$, we have $|U|=|V|\le 4$. Moreover, as $G$ is connected, it follows that
\(
2\le |S|\le 3.
\)

\noindent\textbf{Case 1:} $|S|=2$.

In this case, $|N(S)| = 1$. Since $G$ is connected, the unique vertex in $N(S)$ is
adjacent to both vertices of $S$, which are therefore pendant vertices. Hence,
$G$ contains a vertex adjacent to more than one pendant vertex. By
\thmref{antenna}, $G$ does not have property (P), a contradiction.

\noindent\textbf{Case 2:} $|S|=3$.

Since $S \subseteq U$ and $|V(G)|\le 8$, we have $3 \le |U| \le 4$. If $|U|=3$, then $S=U$. Since $G$ is connected, $N(S)=V$, giving $|N(S)|=|V|=3=|S|$, a contradiction. Hence, $|U|=|V|=4$.

Thus $|N(S)|\le 2$. 
\begin{enumerate}[label=\textup{(\roman*)}]
\item If $|N(S)| = 1$, then the unique vertex in $N(S)$ is adjacent to all three
vertices of $S$, which are pendant vertices. Thus, $G$ is an antenna graph and, by 
\thmref{antenna}, $G$ does not have property (P).

\item If $|N(S)| = 2$, write
\[
U \setminus S = \{u\}, \qquad V \setminus N(S) = \{w_1, w_2\}.
\]
Since $w_1$ and $w_2$ have no neighbors in $S$, the connectivity of $G$ forces $u$ to be
adjacent to both $w_1$ and $w_2$. Moreover, to maintain connectivity, the vertex
$u$ must also be adjacent to at least one vertex in $N(S)$. Hence, $u$ is adjacent
to more than one pendant vertices and, by \thmref{antenna}, $G$ does not have
property (P).
\end{enumerate}
In all cases, the failure of Hall’s condition produces a vertex adjacent to
more than one pendant vertex, contradicting the assumption that $G$ has
property (P). Therefore, $G$ must admit a perfect matching.
\end{pf}
\begin{rem}\label{rem:order8-general}
Let $G$ be a bipartite graph (not necessarily connected) of order at most $8$. Then $G$ has property~(P) if and only if it has a perfect matching.
Indeed, the connected case follows from \lemref{thm:order8}. If $G$ is disconnected, then by \cite[Lemma~2.1]{Howlader}, it suffices to apply \lemref{thm:order8} to each connected component.
\end{rem}
We now extend the characterization obtained for bipartite graphs of order at most $8$ (\lemref{thm:order8}) to a broader class of bipartite graphs
with pendant vertices.
\begin{thm}\label{8+2k}
Let $G$ be a bipartite graph of order $8+2k \ (k \ge 0)$ with at least $k$ pendant vertices. Then $G$ has property (P) if and only if $G$ has a perfect matching.
\end{thm}
\begin{pf}
If $G$ has a perfect matching, then by \thmref{bipartite perfect matching theorem} it follows that $G$ has property (P).

Conversely, suppose that $G$ has property (P). We prove that $G$ has a perfect matching by induction on $k$.

The case $k=0$ follows directly from \remref{rem:order8-general}.

Assume that the statement holds for $k-1$; that is, every bipartite graph of order $8+2(k-1)$ with at least $k-1$ pendant vertices has property (P) if and only if it has a perfect matching.

Now, let $G$ be a bipartite graph of order $8+2k$ with property (P) and containing at least $k$ pendant vertices.

\noindent\textbf{Case 1:} Some vertex of $G$ is adjacent to two or more pendant vertices. In this case, $G$ is an antenna graph. By \thmref{antenna}, $G$ does not have property (P). Hence, this case cannot occur.

\noindent\textbf{Case 2:} Every vertex of $G$ is adjacent to at most one pendant vertex. Since $G$ has at least $k$ pendant vertices, choose an edge $vu$, where $u$ is the pendant vertex and $v$ is its neighbor. 
By \thmref{lem:pendant_heredity}, since $G$ has property (P), $G - \{u,v\}$ has property (P). Moreover,
\[
|V(G - \{u,v\})| = (8+2k) - 2 = 8+2(k-1).
\]
The graph $G - \{u,v\}$ is a bipartite graph with at least $k-1$ pendant vertices and has property (P). By the induction hypothesis, $G - \{u,v\}$ has a perfect matching, say $M'$.
The set of edges
\[
M = M' \cup \{vu\}
\]
is a perfect matching of $G$.
\end{pf}
\begin{figure}[h]
\centering
\begin{tikzpicture}[scale=1]

% vertex style
\tikzset{
  vtx/.style={circle, fill=black, inner sep=1.6pt}
}

%---------------- Left part U ----------------%
\node[vtx] (u1) at (0,5) {}; \node[left=6pt] at (u1) {$u_1$};
\node[vtx] (u2) at (0,4) {}; \node[left=6pt] at (u2) {$u_2$};
\node[vtx] (u3) at (0,3) {}; \node[left=6pt] at (u3) {$u_3$};
\node[vtx] (u4) at (0,2) {}; \node[left=6pt] at (u4) {$u_4$};
\node[vtx] (u5) at (0,1) {}; \node[left=6pt] at (u5) {$u_5$};
\node[vtx] (u6) at (0,0) {}; \node[left=6pt] at (u6) {$u_6$};

%---------------- Right part V ----------------%
\node[vtx] (v1) at (3,5) {}; \node[right=6pt] at (v1) {$v_1$};
\node[vtx] (v2) at (3,4) {}; \node[right=6pt] at (v2) {$v_2$};
\node[vtx] (v3) at (3,3) {}; \node[right=6pt] at (v3) {$v_3$};
\node[vtx] (v4) at (3,2) {}; \node[right=6pt] at (v4) {$v_4$};
\node[vtx] (v5) at (3,1) {}; \node[right=6pt] at (v5) {$v_5$};
\node[vtx] (v6) at (3,0) {}; \node[right=6pt] at (v6) {$v_6$};

%---------------- Edges ----------------%
% K_{3,2} violator block
\draw (u1)--(v1); \draw (u1)--(v2);
\draw (u2)--(v1); \draw (u2)--(v2);
\draw (u3)--(v1); \draw (u3)--(v2);

% K_{2,3} lower block
\draw (u4)--(v3); \draw (u4)--(v4); \draw (u4)--(v5);
\draw (u5)--(v3); \draw (u5)--(v4); \draw (u5)--(v5);

% Bridge edge to make it connected
\draw (u4)--(v2);

\draw (u6)--(v5);
\draw (v6)--(u5);

\end{tikzpicture}
\caption{A connected bipartite graph of order $12$ with cycle rank $4$ and two pendant vertices.}
\label{fig:pendant_example}
\end{figure}

\begin{example}\label{ex:power_of_8_plus_2k}
Consider the graph $G$ in Figure~\ref{fig:pendant_example}, which has order $12$ and cycle rank $4$. Because $G$ has exactly $2$ pendant vertices, $u_6$ and $v_6$, it satisfies the hypotheses of \thmref{8+2k} with $k=2$. Observe that the subset $\{u_1, u_2, u_3\}$ has only two neighbors $\{v_1, v_2\}$, which violates Hall's condition. Therefore, $G$ has no perfect matching. By \thmref{8+2k}, we immediately conclude that $G$ does not have property (P).
\end{example}
Since deleting a pendant vertex and its unique neighbor preserves both perfect matchings and property~(P) (\thmref{lem:pendant_heredity}), we introduce the following reduction process.
\begin{defn}
Let $G$ be a graph. A graph $\bar{G}$ is said to be a
P-reduced graph of $G$ if it is obtained from $G$
by the pendant-vertex reduction process described below.

In each step, select a pendant vertex of the graph $G$ (if one exists) and
delete it together with its unique neighbor. Repeat this process until the resulting graph satisfies one of the following conditions:
\begin{enumerate}[label=\textup{(\roman*)}]
\item it has no pendant vertices;
\item it contains an antenna, i.e., a vertex adjacent to more than one pendant vertex.
\end{enumerate}
Note that the P-reduced graph need not be unique (see the example below), since different choices of pendant vertices during the reduction may lead to different resulting graphs, and any graph obtained in this manner is called a P-reduced graph of $G$.
\end{defn}
\begin{figure}[h]
\centering
\begin{tikzpicture}[scale=1]

\tikzset{
vtx/.style={circle, fill=black, inner sep=1.8pt}
}

%---------------- G ----------------%

\node[vtx,label=left:$u_1$] (u1) at (0,6) {};
\node[vtx,label=left:$u_2$] (u2) at (0,5) {};
\node[vtx,label=left:$u_3$] (u3) at (0,4) {};
\node[vtx,label=left:$u_4$] (u4) at (0,3) {};
\node[vtx,label=left:$u_5$] (u5) at (0,2) {};
\node[vtx,label=left:$u_6$] (u6) at (0,1) {};

\node[vtx,label=right:$v_1$] (v1) at (2,6) {};
\node[vtx,label=right:$v_2$] (v2) at (2,5) {};
\node[vtx,label=right:$v_3$] (v3) at (2,4) {};
\node[vtx,label=right:$v_4$] (v4) at (2,3) {};
\node[vtx,label=right:$v_5$] (v5) at (2,2) {};
\node[vtx,label=right:$v_6$] (v6) at (2,1) {};

\draw (u1)--(v1);
\draw (v1)--(u2);
\draw (v1)--(u4);
\draw (u2)--(v2);
\draw (v2)--(u3);
\draw (u3)--(v3);
\draw (v3)--(u4);

\draw (u5)--(v4);
\draw (u5)--(v5);
\draw (u5)--(v6);
\draw (u6)--(v5);
\draw (u6)--(v4);
\draw (u6)--(v6);
\draw (v3)--(u5);
\draw (v3)--(u6);

\node at (1,0.2) {$G$};

% arrow
\draw[->,thick] (1,-0.2)--(1,-1)
node[right,midway]{\small delete $u_1,v_1$};

%---------------- G' ----------------%

\node[vtx,label=left:$u_2$] (a1) at (0,-2) {};
\node[vtx,label=left:$u_3$] (a2) at (0,-3) {};
\node[vtx,label=left:$u_4$] (a3) at (0,-4) {};
\node[vtx,label=left:$u_5$] (a4) at (0,-5) {};
\node[vtx,label=left:$u_6$] (a5) at (0,-6) {};

\node[vtx,label=right:$v_2$] (b1) at (2,-2) {};
\node[vtx,label=right:$v_3$] (b2) at (2,-3) {};
\node[vtx,label=right:$v_4$] (b3) at (2,-4) {};
\node[vtx,label=right:$v_5$] (b4) at (2,-5) {};
\node[vtx,label=right:$v_6$] (b5) at (2,-6) {};

\draw (a1)--(b1);
\draw (b1)--(a2);
\draw (a2)--(b2);
\draw (b2)--(a3);

\draw (a4)--(b3);
\draw (b2)--(a4);
\draw (a5)--(b2);
\draw (a5)--(b3);
\draw (a5)--(b4);
\draw (a5)--(b5);
\draw (a4)--(b5);
\draw (a4)--(b4);

\node at (1,-7) {$G'$};

% arrows to reduced graphs
\draw[->,thick] (0.9,-7.5)--(-1,-8.7)
node[midway, left=7pt]{\small delete $u_2,v_2$};

\draw[->,thick] (1.1,-7.5)--(3,-8.7)
node[midway, right=7pt]{\small delete $u_4,v_3$};

%---------------- G1 ----------------%

\node[vtx,label=left:$u_3$] (c1) at (-2,-9) {};
\node[vtx,label=left:$u_4$] (c2) at (-2,-10) {};
\node[vtx,label=left:$u_5$] (c3) at (-2,-11) {};
\node[vtx,label=left:$u_6$] (c4) at (-2,-12) {};

\node[vtx,label=right:$v_3$] (d1) at (0,-9) {};
\node[vtx,label=right:$v_4$] (d2) at (0,-10) {};
\node[vtx,label=right:$v_5$] (d3) at (0,-11) {};
\node[vtx,label=right:$v_6$] (d4) at (0,-12) {};

\draw (c1)--(d1);
\draw (d1)--(c2);
\draw (c4)--(d2);
\draw (d2)--(c3);
\draw (c3)--(d3);
\draw (d3)--(c4);
\draw (c4)--(d4);
\draw (c3)--(d1);
\draw (c4)--(d1);
\draw (c3)--(d4);

\node at (-1,-12.7) {$G_1$};

%---------------- G2 ----------------%

\node[vtx,label=left:$u_2$] (e1) at (2,-9) {};
\node[vtx,label=left:$u_3$] (e2) at (2,-10) {};
\node[vtx,label=left:$u_5$] (e3) at (2,-11) {};
\node[vtx,label=left:$u_6$] (e4) at (2,-12) {};

\node[vtx,label=right:$v_2$] (f1) at (4,-9) {};
\node[vtx,label=right:$v_4$] (f2) at (4,-10) {};
\node[vtx,label=right:$v_5$] (f3) at (4,-11) {};
\node[vtx,label=right:$v_6$] (f4) at (4,-12) {};

\draw (e1)--(f1);
\draw (f1)--(e2);
\draw (e3)--(f2);
\draw (f3)--(e3);
\draw (e3)--(f4);
\draw (f3)--(e4);
\draw (e4)--(f2);
\draw (e4)--(f4);

\node at (3,-12.7) {$G_2$};

\end{tikzpicture}

\caption{Pendant-vertex reduction of $G$ yields the P-reduced graphs $G_1$ and $G_2$.}

\label{fig:pendant12}
\end{figure}
\begin{example}\label{ex:pendant12}
Consider the connected bipartite graph $G$ of order $12$ shown in
Figure~\ref{fig:pendant12}. The graph $G$ contains a pendant vertex
$u_1$. Deleting $u_1$ together with its neighbor
$v_1$ yields the graph $G'$. In the graph $G'$, there are two pendant vertices, namely $u_2$ and $u_4$. Deleting $u_2$ together with its neighbor $v_2$ produces the graph $G_1$,
while deleting $u_4$ together with its neighbor $v_3$ produces the graph $G_2$.
Both $G_1$ and $G_2$ contain an antenna. Therefore, the
reduction process terminates at this stage. Consequently, $G_1$ and $G_2$ are P-reduced graphs of $G$ obtained by the pendant-vertex reduction process, so either graph may be denoted
by $\bar{G}$. This example illustrates that the P-reduced graph of
$G$ need not be unique.
\end{example}
The following theorem shows how the pendant-vertex reduction process extends the characterization of property~(P).
\begin{thm}\label{pendant--edge reduction theorem}
Let $G$ be a bipartite graph, and let $\bar{G}$ be a P-reduced graph of $G$
obtained by the pendant-vertex reduction process. If $\bar{G}$ has order at most $8$,
then $G$ has property (P) if and only if $G$ has a perfect matching.
\end{thm}

\begin{pf}
If $G$ has a perfect matching, then by \thmref{bipartite perfect matching theorem} it follows that $G$ has property (P).

Conversely, suppose that $G$ has property (P). We show that $G$ has a perfect matching. If $G$ has no pendant vertices, then by definition $\bar{G}=G$. Since $\bar{G}$ has order at most $8$, \remref{rem:order8-general} implies that $G$ has a perfect matching. If $G$ contains an antenna, then $G=\bar{G}$ and by \thmref{antenna} the graph $G$ does not have property (P), contradicting our assumption. Hence, this case cannot occur.

Now suppose that $G$ has at least one pendant vertex. Repeatedly delete a
pendant vertex together with its unique neighbor. By repeated applications of
\thmref{lem:pendant_heredity}, every graph obtained during this process has
property~(P). Consequently, none of these graphs can contain an antenna, by
\thmref{antenna}. Since each step removes exactly two vertices, the process
terminates after finitely many steps at a graph $\bar G$ with no pendant
vertices.

By hypothesis, $|V(\bar G)|\le8$. Hence, by
\remref{rem:order8-general}, $\bar G$ admits a perfect matching, say
$M_{\bar G}$. Let $e_1,e_2,\ldots,e_t$ be the edges removed during the
reduction process. Since each $e_i$ joins a pendant vertex to its unique
neighbor, the edges $e_1,e_2,\ldots,e_t$ are pairwise vertex-disjoint and
therefore form a matching. Thus
\[
M=M_{\bar G}\cup\{e_1,e_2,\ldots,e_t\}
\]
is a perfect matching of $G$.
\end{pf}
\begin{example}
Consider the graph $G$ shown in Figure~\ref{fig:pendant12}, which has cycle rank $4$. In Example~\ref{ex:pendant12},
the pendant-vertex reduction process yields P-reduced graphs $G_1$
and $G_2$, each of which has order $8$. Moreover, $G$ does not have a perfect matching, since the vertex set
$\{v_4,v_5,v_6\}$ violates Hall's condition. Consequently, by
Theorem~\ref{pendant--edge reduction theorem}, the graph $G$ does not
have property (P).
\end{example}
Note that the preceding example cannot be addressed using Theorem~\ref{8+2k}. Since the graph has order $12$ but contains only one pendant vertex. This explicitly highlights the broader applicability of Theorem~\ref{pendant--edge reduction theorem}.
%%%%%%%%%%%%%%%%%%%%%%%%%%%%%%%%%%%%%%%%%%%%%%%%%%%%%%%%%%%%%%%%%%%%%%%%%%%%%%%%%%%%%%%%%%%%%%%%%%%%%%%%%%%%%%%%
%%%%%%%%%%%%%%%%%%%%%%%%%%%%%%%%%%%%%%%%%%%%%%%%%%%%%%%%%%%%%%%%%%%%%%%%%%%%%%%%%%%%%%%%%%%%%%%%%%%%%%%%%%%%%%%%
\section{Triangular Bipartite Graphs}\label{sec:triangular}
We study bipartite graphs whose biadjacency matrices are in upper triangular form and characterize when such graphs have property (P).

\begin{defn}[Triangular bipartite graph]\label{def:triangular}
Let $G=G(U,V;E)$ be a balanced bipartite graph with $|U|=|V|=n$. 
We call $G$ \emph{triangular} if there exist orderings
\[
U=\{u_1,\dots,u_n\} \ \text{ and } \ V=\{v_1,\dots,v_n\},
\]
such that its biadjacency matrix $B$ is upper triangular, i.e.,  it satisfies
\(B_{ij}=0 \ \text{for all } i>j.
\)
\end{defn}
\begin{rem}
Because the biadjacency matrix is upper triangular, every perfect matching in $G$ must correspond exactly to the main diagonal of the matrix. Consequently, a triangular bipartite graph admits at most one perfect matching.
\end{rem}
\begin{thm}\label{triangular-theorem}
A triangular bipartite graph has property (P) if and only if it admits a perfect matching.
\end{thm}
\begin{pf}
If $G$ has a perfect matching, then by \thmref{bipartite perfect matching theorem}, $G$ has property (P).

Conversely, suppose that $G$ has property~(P). Let $U$ and $V$ be ordered as in Definition~\ref{def:triangular}. Then $G$ admits a hollow-inverse matrix
\[
\renewcommand{\arraystretch}{1.3}
A=
\left[
\begin{array}{c|c}
X & D\\ \hline
D^{\top} & Y
\end{array}
\right],
\]
where $X$ and $Y$ are diagonal matrices, and $D$ has the zero-nonzero pattern defined by $G$. Since $G$ is triangular, $D$ is upper triangular.

Write
\[
\renewcommand{\arraystretch}{1.3}
A^{-1}=
\left[
\begin{array}{c|c}
P & Q\\ \hline
Q^{\top} & R
\end{array}
\right].
\]
Since $A^{-1}$ is hollow,
\[
P_{ii}=R_{ii}=0, \ \text{for all } 1\le i\le n.
\]

From $AA^{-1}=I_{2n}$, we obtain
\begin{equation}\label{1st diag}
    XP + DQ^{\top} = I_n,
\end{equation}
\begin{equation}\label{2nd diag}
    D^{\top}Q +YR = I_n,
\end{equation}
\begin{equation}\label{1st off}
    XQ + DR = \mathbf{0}_{n},
\end{equation}
\begin{equation}\label{2nd off}
    D^{\top}P + YQ^{\top}= \mathbf{0}_{n},
\end{equation}
where $\mathbf{0}_n$ denotes the zero matrix of order $n$.

We prove by induction that
\[
D_{ii}\neq0,\qquad 1\le i\le n.
\]

For $i=1$, the $(1,1)$-entry of \eqref{2nd diag} gives $
(D^{\top}Q)_{11}=1,
$
since $R_{11}=0$ and $Y$ is a diagonal matrix. Since $D$ is upper triangular,
\[
D_{11}Q_{11}=1.
\]
Hence
\[
D_{11}\neq0
\quad\text{and}\quad
Q_{11}\neq0.
\]

From the $(1,1)$-entry of \eqref{2nd off}, using $P_{11}=0$ and $Q_{11}\neq 0$, we obtain $Y_{11}=0$. Then the first row of \eqref{2nd off} and \eqref{2nd diag} yields
\[
P_{1j}=0 \ \text{ for } 1\le j\le n, \qquad Q_{1j}=0 \ \text{ for } 1<j\le n.
\] 

Assume that, for some $1\le k<n$,
\[
D_{ii}\neq0 \ \text{for all } i\le k,
\]
and
\[
P_{ij}=0 \ \text{for all } i\le k,\ 1\le j\le n,
\quad 
Q_{ij}=0 \ \text{for all } i\le k,\ i<j\le n.
\]
Then the $(k+1,k+1)$-entry of \eqref{2nd diag} gives
\[
D_{(k+1)(k+1)}Q_{(k+1)(k+1)}=1,
\]
since $Q_{i,k+1}=0$ for $i\le k$ and $R_{(k+1)(k+1)}=0$. 

Hence, \[
D_{(k+1)(k+1)}\neq0,
\qquad
Q_{(k+1)(k+1)}\neq0.
\]

Taking the $(k+1,k+1)$-entry of \eqref{2nd off}, it follows that $Y_{(k+1)(k+1)}=0$. Then the $(k+1)$-th row of \eqref{2nd off} and \eqref{2nd diag} yields
\[
P_{(k+1)j}=0 \ \text{for all } 1\le j\le n, \qquad Q_{(k+1)j}=0 \ \text{for } k+1<j\le n.
\]
This completes the induction.

Thus $D_{ii}\neq0$ for every $i$, so
\[
M=\{u_1v_1,\ldots,u_nv_n\}
\]
is a perfect matching of $G$.
\end{pf}
\begin{figure}[h]
\centering
\begin{tikzpicture}[scale=1]

% vertex style
\tikzset{
  vtx/.style={circle, fill=black, inner sep=1.6pt}
}

%---------------- Left part U ----------------%
\node[vtx] (u1) at (0,7) {}; \node[left=6pt] at (u1) {$u_1$};
\node[vtx] (u2) at (0,6) {}; \node[left=6pt] at (u2) {$u_2$};
\node[vtx] (u3) at (0,5) {}; \node[left=6pt] at (u3) {$u_3$};
\node[vtx] (u4) at (0,4) {}; \node[left=6pt] at (u4) {$u_4$};
\node[vtx] (u5) at (0,3) {}; \node[left=6pt] at (u5) {$u_5$};
\node[vtx] (u6) at (0,2) {}; \node[left=6pt] at (u6) {$u_6$};
\node[vtx] (u7) at (0,1) {}; \node[left=6pt] at (u7) {$u_7$};
\node[vtx] (u8) at (0,0) {}; \node[left=6pt] at (u8) {$u_8$};

%---------------- Right part V ----------------%
\node[vtx] (v1) at (3,7) {}; \node[right=6pt] at (v1) {$v_1$};
\node[vtx] (v2) at (3,6) {}; \node[right=6pt] at (v2) {$v_2$};
\node[vtx] (v3) at (3,5) {}; \node[right=6pt] at (v3) {$v_3$};
\node[vtx] (v4) at (3,4) {}; \node[right=6pt] at (v4) {$v_4$};
\node[vtx] (v5) at (3,3) {}; \node[right=6pt] at (v5) {$v_5$};
\node[vtx] (v6) at (3,2) {}; \node[right=6pt] at (v6) {$v_6$};
\node[vtx] (v7) at (3,1) {}; \node[right=6pt] at (v7) {$v_7$};
\node[vtx] (v8) at (3,0) {}; \node[right=6pt] at (v8) {$v_8$};

%---------------- Edges (upper triangular) ----------------%

\foreach \i/\j in {
1/1,1/2,1/3,1/4,1/5,
2/2,2/7,2/8,
3/3, 3/4,3/6,3/8,
4/5,4/6,4/7,4/8,
5/6,5/8,
6/7,6/8,
7/7,7/8,
8/8}
{
  \draw (u\i)--(v\j);
}

\end{tikzpicture}
\caption{A triangular bipartite graph of order $16$.}
\label{fig:triangular16}
\end{figure}

\begin{example}
Figure~\ref{fig:triangular16} shows a triangular bipartite graph $G$ of order $16$ and cycle rank $8$, with ordered partite sets
\(
U=\{u_1,\dots,u_8\}, \ \text{and } V=\{v_1,\dots,v_8\},
\)
for which the biadjacency matrix $B$ is upper triangular and is given by
\[
B=
\begin{bmatrix}
1 & 1 & 1 & 1 & 1 & 0 & 0 & 0\\
0 & 1 & 0 & 0 & 0 & 0 & 1 & 1\\
0 & 0 & 1 & 1 & 0 & 1 & 0 & 1\\
0 & 0 & 0 & 0 & 1 & 1 & 1 & 1\\
0 & 0 & 0 & 0 & 0 & 1 & 0 & 1\\
0 & 0 & 0 & 0 & 0 & 0 & 1 & 1\\
0 & 0 & 0 & 0 & 0 & 0 & 1 & 1\\
0 & 0 & 0 & 0 & 0 & 0 & 0 & 1
\end{bmatrix}.
\]
Hence $G$ is a triangular bipartite graph.
However, $G$ does not admit a perfect matching.
Consider the subset
\(
S=\{u_6,u_7,u_8\}\subseteq U.
\) Then 
its set of neighbors is 
\(
N(S)=\{v_7,v_8\},
\) and $|N(S)|=2<3=|S|$.
Thus, Hall's condition fails, and $G$ has no perfect matching. Now, by Theorem~\ref{triangular-theorem},
$G$ does not have property (P).
\end{example}
%%%%%%%%%%%%%%%%%%%%%%
\section{Matrix Constraints on Bipartite Graphs Without Perfect Matchings}
\label{sec:matrix_constraints}

In this final section, we derive necessary algebraic conditions for hollow-inverse matrices on bipartite graphs without perfect
matchings and use them to identify a class of bipartite
graphs that do not have property~(P).

Recall that \(A(S)\) denotes the principal submatrix obtained from \(A\) by deleting the rows and columns indexed by \(S\). We also write \(A[S]:=A[S,S]\) for the principal submatrix induced by \(S\), and \(A[S,T]\) for the submatrix with rows indexed by \(S\) and columns indexed by \(T\).
%%%%%%%%%%%%%%%%%%%%%%%%%%%%%%%%%%%%%%%%%%%%%%%%%%%%%%%%%%%%%%%%%%%
\subsection{General Algebraic Constraints for Bipartite Graphs Lacking Perfect Matchings}
Let \(G=G(U,V;E)\) be a balanced bipartite graph of order \(2n\) with no
perfect matching. By Hall's theorem, there exists a subset
\(U_1\subseteq U\) such that
\[
 |N(U_1)|<|U_1|.
\]
Let \(V_1=N(U_1)\), and write
\[
 |U_1|=r, \quad |V_1|=k, \quad r>k.
\]
Set
\[
 U_2=U\setminus U_1,\qquad V_2=V\setminus V_1,
 \qquad s=|V_2|=n-k,
\]
and denote the Hall deficiency by
\[
 d=r-k\ge 1.
\]
Since \(U_1\) has no neighbors in \(V_2\), the rows and columns of any
matrix \(A\in S(G)\), ordered according to the partition
\((U_1,U_2,V_1,V_2)\), have the block form
\[
\renewcommand{\arraystretch}{1.3}
A=
\left[
\begin{array}{cc|cc}
X_1 & \mathbf{0} & D_{11} & \mathbf{0}\\
\mathbf{0} & X_2 & D_{21} & D_{22}\\ \hline
D_{11}^{\top} & D_{21}^{\top} & Y_1 & \mathbf{0}\\
\mathbf{0} & D_{22}^{\top} & \mathbf{0} & Y_2
\end{array}
\right],
\]
where
\[
D_{11}=A[U_1,V_1],\qquad
D_{21}=A[U_2,V_1],\qquad
D_{22}=A[U_2,V_2].
\]
Here \(X_i\) and \(Y_i\) are diagonal matrices containing the diagonal
entries of \(A\) corresponding to the respective vertex subsets.

\begin{lem}\label{lem:weights}
Let \(A\in S(G)\) be a hollow-inverse matrix, and let \(C=A^{-1}\).
For each edge \(ij\in E(G)\), define
\[
\omega_{ij}:=A_{ij}C_{ij}.
\]
Then
\begin{equation}\label{eq:vertexsum}
 \sum_{j:\,ij\in E(G)} \omega_{ij}=1
 \qquad\text{for every vertex } i.
\end{equation}
\end{lem}

\begin{pf}
Since \(AC=I\), for each vertex \(i\) we have
\[
1=(AC)_{ii}
 = \sum_j A_{ij}C_{ji}
 = A_{ii}C_{ii}+\sum_{j\ne i} A_{ij}C_{ij}.
\]
The first term is zero because \(C\) is hollow. Moreover,
\(A_{ij}=0\) whenever \(ij\notin E(G)\). Hence
\[
1=\sum_{j:\,ij\in E(G)} A_{ij}C_{ij}
 =\sum_{j:\,ij\in E(G)} \omega_{ij},
\]
as required.
\end{pf}

\begin{prop}\label{prop:cut}
Suppose \(G\) has no perfect matching, and let \(U_1\subseteq U\) be a
subset violating Hall's condition. Let \(A\in S(G)\) be a hollow-inverse
matrix with inverse \(C=A^{-1}\). For subsets \(X,Y\subseteq V(G)\), define
the cross-weight
\[
\omega(X,Y)
:=
\sum_{\substack{u\in X,\;w\in Y\\ uw\in E(G)}} A_{uw}C_{uw}.
\]
Then
\[
\omega(U_2,V_1)=k-r\le -1.
\]
In particular, if there is a unique edge \(u^\ast v^\ast\) between 
\(U_2\) and \(V_1\), with \(u^\ast\in U_2\) and \(v^\ast\in V_1\), then
\[
A_{u^\ast v^\ast}C_{u^\ast v^\ast}=-(r-k).
\]
\end{prop}

\begin{pf}
By Lemma~\ref{lem:weights}, the sum of the weights on the edges incident
with each vertex is equal to \(1\).
Summing over all vertices of \(U_1\), and using the fact that every
neighbor of a vertex in \(U_1\) lies in \(V_1\), gives
\[
\omega(U_1,V_1)=r.
\]
On the other hand, summing over all vertices of \(V_1\) yields
\[
\omega(U_1,V_1)+\omega(U_2,V_1)=k.
\]
Therefore,
\[
\omega(U_2,V_1)=k-r\le -1,
\]
since \(r>k\).

If \(u^\ast v^\ast\) is the unique edge between \(U_2\) and \(V_1\), then
\[
\omega(U_2,V_1)=A_{u^\ast v^\ast}C_{u^\ast v^\ast}.
\]
Hence
\[
A_{u^\ast v^\ast}C_{u^\ast v^\ast}=k-r=-(r-k).
\]
\end{pf}
The following proposition provides a lower bound on the number of nonzero diagonal entries of a hollow-inverse matrix when the graph has no perfect matching. This complements a recent result of Sharma and Panda~\cite[Theorem 3.12]{sharma}.
\begin{prop}\label{lem:nondeg}
Let \(A\in S(G)\) be a hollow-inverse matrix with inverse \(C=A^{-1}\). Then
\[
\bigl|\{p\in U_1 : A_{pp}=0\}\bigr|
\le \operatorname{rank}(D_{11})
\le k,
\]
and
\[
\bigl|\{q\in V_2 : A_{qq}=0\}\bigr|
\le \operatorname{rank}(D_{22})
\le n-r.
\]
Consequently, at least \(r-k\) vertices of \(U_1\), and at least \(r-k\)
vertices of \(V_2\), have nonzero diagonal entries.
\end{prop}

\begin{pf}
Let
\[
Z=\{p\in U_1:A_{pp}=0\}.
\]
Fix \(p\in Z\), and let \(c=C e_p\) be the \(p\)-th column of \(C\).
Since \(C\) is hollow, \(C_{pp}=0\), and hence
\[
c=\sum_{j\ne p} C_{jp}e_j.
\]
Using \(AC=I\), we obtain
\[
e_p=A(Ce_p)=Ac=\sum_{j\ne p} C_{jp}Ae_j.
\]
Thus
\[
e_p\in \operatorname{span}\{Ae_j:j\ne p\}.
\]

Now restrict this containment to the rows indexed by \(U_1\). The columns
of \(A\) indexed by \(U_2\) and \(V_2\) vanish on these rows, while the
columns indexed by \(V_1\) span the column space \(\operatorname{Col}(D_{11})\). Moreover,
since \(A[U_1]\) is diagonal, the columns indexed by \(U_1\setminus Z\)
contribute precisely the scaled standard basis vectors
\[
\{A_{jj}e_j:j\in U_1\setminus Z\}.
\]
Since \(A_{jj}\ne0\) for all \(j\in U_1\setminus Z\), these vectors span the subspace 
\(
\operatorname{span}\{e_j:j\in U_1\setminus Z\}.
\)

Therefore, for every \(p\in Z\),
\[
e_p\in
\operatorname{span}\{e_j:j\in U_1\setminus Z\}
+
\operatorname{Col}(D_{11}).
\]
It follows that
\[
\mathbb{R}^{|U_1|}
=
\operatorname{span}\{e_j:j\in U_1\setminus Z\}
+
\operatorname{Col}(D_{11}).
\]
Hence
\[
r\le (r-|Z|)+\operatorname{rank}(D_{11}),
\]
and so
\[
|Z|\le \operatorname{rank}(D_{11})\le k.
\]
Consequently,
\[
\bigl|\{p\in U_1:A_{pp}\ne0\}\bigr|
=r-|Z|\ge r-k.
\]

The argument for \(V_2\) is identical, with \(D_{22}\) replacing
\(D_{11}\). Since
\[
\operatorname{rank}(D_{22})\le |U_2|=n-r,
\]
at least
\[
|V_2|-(n-r)=(n-k)-(n-r)=r-k
\]
vertices of \(V_2\) must have nonzero diagonal entries.
\end{pf}

\begin{rem}
\propref{lem:nondeg} complements a recent theorem of Sharma and Panda
\cite[Theorem 3.12]{sharma}, who proved that a bipartite graph admitting a
hollow-inverse matrix with at most three nonzero diagonal entries must have
a perfect matching.
In contrast, \propref{lem:nondeg} shows that if a bipartite graph has no
perfect matching, then every hollow-inverse matrix on the graph must have at
least
\[
2(r-k)\ge 2
\]
nonzero diagonal entries. Thus the first unresolved case occurs when the
number of nonzero diagonal entries is \(4\).
\end{rem}

\subsection{A Structural Obstruction for a Family of Bipartite Graphs}
\label{specific_example}
We now identify a family of balanced bipartite graphs that do not have
property~(P). These graphs contain a subset that violates Hall's condition with deficiency one and is connected to the rest of the graph by a unique edge.
Using the algebraic constraints established above, we show that this
structure forces a contradiction.

\begin{thm}\label{thm:generalized_obstruction}
Let \(G=G(U,V;E)\) be a balanced bipartite graph. Suppose that there exists a
subset \(U_1\subseteq U\), with \(V_1=N(U_1)\), such that
\begin{enumerate}[label=\textup{(\roman*)}]
    \item \(|V_1|=k\le 2\) and \(|U_1|=k+1\);
    \item there is exactly one edge between \(U_2=U\setminus U_1\) and \(V_1\).
\end{enumerate}
Then \(G\) does not have property (P).
\end{thm}

\begin{pf}
If \(k=1\), then \(|V_1|=1\) and \(|U_1|=2\). In this case, \(G\) contains either
 an isolated vertex or an antenna. If \(G\) has an isolated vertex, then
by Lemma~\ref{lem:isolated_vertex}, it does not have property~(P).
If \(G\) contains an antenna, then Theorem~\ref{antenna} implies that it
does not have property~(P). Hence, the result holds when \(k=1\).

We may therefore assume that \(k=2\). Suppose, for a contradiction, that \(G\) has property~(P). Let
\(A\in S(G)\) be a hollow-inverse matrix with \(C=A^{-1}\). 
Let \(u^\ast v^\ast\) be the unique edge between \(U_2\) and \(V_1\), where
\(u^\ast\in U_2\) and \(v^\ast\in V_1\), and let
\(\Delta=\operatorname{diag}(A)\).

Since \(|U_1|=3\) and \(|V_1|=2\), we write
\[
\Delta[U_1]
=
\operatorname{diag}(x_{u_1},x_{u_2},x_{u_3}),
\qquad
\Delta[V_1]
=
\operatorname{diag}(y_{v_1},y_{v_2}).
\]

We first show that \(C[V_1]=\mathbf0\). Consider the
\((U_1,V_1)\)-block of the equation \(AC=I\). Since \(U_1\) has no
neighbors in \(V_2\), this block equation is
\[
A[U_1]C[U_1,V_1]+A[U_1,V_1]C[V_1]=\mathbf{0}.
\]
Equivalently, with \(D_{11}=A[U_1,V_1]\),
\[
\Delta[U_1]C[U_1,V_1]=-D_{11}C[V_1].
\]
Multiplying on the left by \(C[U_1,V_1]^{\top}\), we get
\[
C[U_1,V_1]^{\top}\Delta[U_1]C[U_1,V_1]
=
-C[U_1,V_1]^{\top}D_{11}C[V_1].
\]
Since \(\Delta[U_1]\) is diagonal, the left-hand side is symmetric, and hence so is
\(
C[U_1,V_1]^{\top}D_{11}C[V_1].
\)

On the other hand, consider the \((V_1,V_1)\)-block of \(AC=I\). We have
\[
A[V_1,U_1]C[U_1,V_1]
+
A[V_1]C[V_1]
+
A[V_1,U_2]C[U_2,V_1]
=
I_2.
\]
Let
\[
L=D_{11}^{\top}C[U_1,V_1].
\]
Then
\[
L+\Delta[V_1]C[V_1]+D_{21}^{\top}C[U_2,V_1]=I_2.
\]
Since \(D_{21}^{\top}\) has exactly one nonzero entry corresponding to
\(u^\ast v^\ast\), we may write
\[
D_{21}^{\top}
=
\gamma e_{v^\ast}e_{u^\ast}^{\top},
\qquad
\gamma=A_{u^\ast v^\ast}\ne0.
\]
Hence \(D_{21}^{\top}C[U_2,V_1]\) has only one nonzero row,
\(\gamma\,C[\{u^\ast\},V_1]\), corresponding to \(v^\ast\).

Since \(|U_1|-|V_1|=1\), Proposition~\ref{prop:cut} gives
\[
\gamma C_{u^\ast v^\ast}=-1.
\]
Since \(C\) is symmetric and hollow,
\[
C[V_1]=
\begin{bmatrix}
0&c\\
c&0
\end{bmatrix},
\]
for some scalar \(c\).

Without loss of generality, assume that \(v^\ast=v_1\). From
\[
L=I_2-\Delta[V_1]C[V_1]-D_{21}^{\top}C[U_2,V_1],
\]
we obtain
\[
L
=
\begin{bmatrix}
2 & -cy_{v_1}-\gamma C_{u^\ast v_2}\\
-cy_{v_2} & 1
\end{bmatrix}.
\]
Hence
\[
L^{\top}
=
\begin{bmatrix}
2 & -cy_{v_2}\\
-cy_{v_1}-\gamma C_{u^\ast v_2} & 1
\end{bmatrix}.
\]
Therefore
\[
L^{\top}C[V_1]
=
\begin{bmatrix}
-c^2y_{v_2} & 2c\\
c & c(-cy_{v_1}-\gamma C_{u^\ast v_2})
\end{bmatrix}.
\]
Since
\[
L^{\top}C[V_1]
=
C[U_1,V_1]^{\top}D_{11}C[V_1]
\]
is symmetric, comparing the off-diagonal
entries gives
\(
2c=c.
\)
Hence \(c=0\), so \(C[V_1]=\mathbf{0}\).

With \(C[V_1]=\mathbf{0}\), the \((U_1,V_1)\)-block equation reduces to
\[
\Delta[U_1]C[U_1,V_1]=\mathbf{0}.
\]
Moreover,
\[
L
=
D_{11}^{\top}C[U_1,V_1]
=
\begin{bmatrix}
2 & -\gamma C_{u^\ast v_2}\\
0 & 1
\end{bmatrix}.
\]
In particular, \(L\) is invertible. Hence \(C[U_1,V_1]\) has full column
rank \(2\). Since
\[
\Delta[U_1]C[U_1,V_1]=\mathbf{0},
\]
the nullity of the diagonal matrix \(\Delta[U_1]\) must be at least \(2\).
Thus at least two of the three diagonal entries of \(\Delta[U_1]\) are zero.

Now consider the \((U_1,U_1)\)-block of \(AC=I\),
\[
\Delta[U_1]C[U_1]+D_{11}C[V_1,U_1]=I_3.
\]
Equivalently,
\[
I_3-\Delta[U_1]C[U_1]=D_{11}C[V_1,U_1].
\]
The right-hand side is a product of a \(3\times2\) matrix and a
\(2\times3\) matrix, and therefore has rank at most \(2\).

By Proposition~\ref{lem:nondeg}, \(\Delta[U_1]\) has at least one nonzero diagonal entry. Since at least two of its diagonal entries are zero, exactly one diagonal entry is nonzero. Without loss of generality, assume that
\(
x_{u_3}\ne0.
\)

Because \(C\) is hollow,
\[
C[U_1]=
\begin{bmatrix}
0 & C_{12} & C_{13}\\
C_{12} & 0 & C_{23}\\
C_{13} & C_{23} & 0
\end{bmatrix}.
\]
Hence
\[
\Delta[U_1]C[U_1]
=
\begin{bmatrix}
0 & 0 & 0\\
0 & 0 & 0\\
x_{u_3}C_{13} & x_{u_3}C_{23} & 0
\end{bmatrix}.
\]
This matrix is nilpotent. Therefore every eigenvalue of
\(
I_3-\Delta[U_1]C[U_1]
\)
is equal to \(1\), and hence
\[
\operatorname{rank}\bigl(I_3-\Delta[U_1]C[U_1]\bigr)=3.
\]
This is impossible, since
\[
I_3-\Delta[U_1]C[U_1]
=
D_{11}C[V_1,U_1],
\]
where the right-hand side has rank at most \(2\). Therefore no such
hollow-inverse matrix \(A\) exists, and \(G\) does not have
property~(P).
\end{pf}
We now return to the graph from \remref{rem:sharp} (\secref{sec:schur_reductions}) and show that it does not have property~(P).

\begin{cor}\label{cor:rank_contradiction}
The graph \(G\) in Figure~\ref{fig:sharp-example} does not have property~(P).
\end{cor}
\begin{pf}
Let
\(
U_1=\{u_1,u_2,u_3\},
\ 
V_1=N(U_1)=\{v_1,v_2\}.
\)
Then
\(
|U_1|=3,\  |V_1|=2,
\)
and there is exactly one edge, namely \(u_4v_1\), between
\(U_2=U\setminus U_1\) and \(V_1\). Therefore \(G\) satisfies the hypotheses of
Theorem~\ref{thm:generalized_obstruction} with \(k=2\), and hence does not
have property~(P).
\end{pf}

\begin{rem}\label{rem:limitations}
The assumptions \(k\le2\) and the existence of a unique edge between
\(U_2\) and \(V_1\) are essential in the proof of
Theorem~\ref{thm:generalized_obstruction}. If more than one edge lies
between \(U_2\) and \(V_1\), the symmetry argument no longer yields the key
identity used in the proof. Likewise, when \(k\ge3\), a hollow symmetric
matrix may contain a zero column without being the zero matrix. Consequently,
the argument no longer forces \(C[V_1]=\mathbf{0}\), and the subsequent rank
contradiction fails.
It remains an open problem whether similar obstructions can be established
when \(k\ge3\) or when there is more than one edge between \(U_2\) and
\(V_1\).
\end{rem}
%%%%%%%%%%%%%%%%%%%%%%

\section*{Declaration of competing interest}
There is no competing interest.

\section*{Funding}
The first author's research was supported by the Anusandhan National Research Foundation (ANRF), Government of India, under Startup Research Grant SRG/2022/001281 and MATRICS Grant ANRF/ARGM/2025/002777/MTR. The second author's research was supported by the National Board for Higher Mathematics (NBHM), Government of India, under the Ph.D. Scholarship Grant 0203/8(28)/2023-R\&D-II.
%%%%%%%%%%%%%%%%%%%%%%%%%%%%%%%%%%%%%%%%%%%%%%%%%%%%%%%%%%%%%%%%%%%%%%%%%%%%%%%%%%%%%%%%%%%%%%%%%%%%%%%%%%%%


\begin{thebibliography}{99}
\bibitem{fons1}M. Andeli\'c, C.~M. da~Fonseca and R. Mamede, On the number of P-vertices of some graphs, Linear Algebra Appl. {\bf 434} (2011), no.~2, 514--525; MR2741238

\bibitem{fons2}M. Andeli\'c, A.~L. Eri\'c{} and C.~M. da~Fonseca, Nonsingular acyclic matrices with full number of P-vertices, Linear Multilinear Algebra {\bf 61} (2013), no.~1, 49--57; MR3003041

\bibitem{cvetkovic2010}Cvetković, Dragoš, Peter Rowlinson, and Slobodan Simić. "An introduction to the theory of graph spectra." London Mathematical Society Student Texts 75 (2010).

\bibitem{diestel}R. Diestel, {\it Graph theory}, sixth edition, 
Graduate Texts in Mathematics, 173, Springer, Berlin, [2025] \copyright 2025; MR4874150

\bibitem{fons5} Z. Du and C.~M. da~Fonseca, Nonsingular acyclic matrices with an extremal number of P-vertices, Linear Algebra Appl. {\bf 442} (2014), 2--19; MR3134347

\bibitem{fons3}Z. Du and C.~M. da~Fonseca, The singular acyclic matrices with maximal number of P-vertices, Linear Algebra Appl. {\bf 438} (2013), no.~5, 2274--2279; MR3005289

\bibitem{du1} Z. Du and C.~M. da~Fonseca, The acyclic matrices with a P-set of maximum size, Linear Algebra Appl. {\bf 468} (2015), 27--37; MR3293238

\bibitem{du2} Z. Du and C.~M. da~Fonseca, The number of P-vertices for acyclic matrices with given nullity, Discrete Math. {\bf 346} (2023), no.~12, Paper No. 113592, 17 pp.; MR4615398

\bibitem{du5} Z. Du and C.~M. da~Fonseca, The number of P-vertices for acyclic matrices of maximum nullity, Discrete Appl. Math. {\bf 269} (2019), 211--219; MR4016599

\bibitem{du3} Z. Du and C.~M. da~Fonseca, The real symmetric matrices of odd order with a P-set of maximum size, Czechoslovak Math. J. {\bf 66(141)} (2016), no.~3, 1007--1026; MR3556881

\bibitem{du4} Z. Du and C.~M. da~Fonseca, The real symmetric matrices with a given rank and a P-set with maximum size, Discrete Math. {\bf 348} (2025), no.~11, Paper No. 114572, 8 pp.; MR4910416

\bibitem{fons4} A.~L. Eri\'c{} and C.~M. da~Fonseca, The maximum number of P-vertices of some nonsingular double star matrices, Discrete Math. {\bf 313} (2013), no.~20, 2192--2194; MR3084262

\bibitem{cruz} M.~R. Fernandes and H.~F. da~Cruz, The number of $P$-vertices in a matrix with maximum nullity, Linear Algebra Appl. {\bf 547} (2018), 168--182; MR3781367

\bibitem{fons6} A. Fonseca, Â. Mestre, A. Mohammadian, C. Perdigão and M.~M. Torres, The maximum number of Parter vertices of acyclic matrices, Discrete Math. {\bf 344} (2021), no.~2, 112198.

\bibitem{grid perfect matching reference}G.~D. da~Fonseca, B. Ries and D. Sasaki, On the ratio between maximum weight perfect matchings and maximum weight matchings in grids, Discrete Appl. Math. {\bf 207} (2016), 45--55; MR3497983

\bibitem{godsil}C.~D. Godsil, Algebraic matching theory, Electron. J. Combin. {\bf 2} (1995), Research Paper 8, approx.\ 14 pp.; MR1323978

\bibitem{harary}F. Harary, {\it Graph theory}, Addison-Wesley Publishing Co., Reading, Mass.-Menlo Park, Calif.-London, 1969; MR0256911

\bibitem{inter} R.~A. Horn and C.~R. Johnson, {\it Matrix analysis}, second edition, 
Cambridge Univ. Press, Cambridge, 2013; MR2978290

\bibitem{Howlader} A. Howlader, P.~R. Raickwade and K.~C. Sivakumar, The full P-vertex problem for unicyclic graphs, Linear Algebra Appl. {\bf 713} (2025), 74--89; MR4880009

\bibitem{Fiedler} C.~R. Johnson, A. Leal-Duarte and C.~M. Saiago, The Parter-Wiener theorem: refinement and generalization, SIAM J. Matrix Anal. Appl. {\bf 25} (2003), no.~2, 352--361; MR2047422

\bibitem{hermitian} C.~R. Johnson and B.~D. Sutton, Hermitian matrices, eigenvalue multiplicities, and eigenvector components, SIAM J. Matrix Anal. Appl. {\bf 26} (2004/05), no.~2, 390--399; MR2124154

\bibitem{kim1} I.-J. Kim and B.~L. Shader, Non-singular acyclic matrices, Linear Multilinear Algebra {\bf 57} (2009), no.~4, 399--407; MR2522851

\bibitem{kim2}I.-J. Kim and B.~L. Shader, On Fiedler- and Parter-vertices of acyclic matrices, Linear Algebra Appl. {\bf 428} (2008), no.~11-12, 2601--2613; MR2416575

\bibitem{Parter} S.~V. Parter, On the eigenvalues and eigenvectors of a class of matrices, J. Soc. Indust. Appl. Math. {\bf 8} (1960), 376--388; MR0112894

\bibitem{sharma} K. Sharma and S.~K. Panda, The P-vertex problem for symmetric matrices whose associated graphs admit perfect matchings, Linear Algebra Appl. {\bf 731} (2026), 109--138; MR4989797

\bibitem{wiener} G.~M. Wiener, Spectral multiplicity and splitting results for a class of qualitative matrices, Linear Algebra Appl. {\bf 61} (1984), 15--29; MR0755246

\bibitem{west}D.~B. West, {\it Introduction to graph theory}, Prentice Hall, Upper Saddle River, NJ, 1996; MR1367739

\end{thebibliography}
\end{document}